\documentstyle[12pt]{article}

\begin{document}
\setlength{\textwidth}{5.5in}
\setlength{\textheight}{8.5in}
\newtheorem{t.}{Theorem}[section]
\newtheorem{d.}[t.]{Definition}
\newtheorem{l.}[t.]{Lemma}
\newtheorem{p.}[t.]{Proposition}
\newtheorem{c.}[t.]{Corollary}
\newtheorem{e.}[t.]{Example}

\title {A Survey on the Complemented Subspace Problem
} 
\author{Mohammad Sal Moslehian\\
Dept. of Math., Ferdowsi Univ., P.O.Box 1159, Mashhad 91775, Iran\\
E-mail: msalm@math.um.ac.ir\\
http://www.um.ac.ir/$\sim$moslehian/}
\date{}
\maketitle

\begin{abstract} The complemented subspace problem asks, in general, which closed subspaces $M$ of a Banach space $X$ are complemented; i.e. there exists a closed subspace $N$ of $X$ such that $X=M\oplus N$? This problem is in the heart of the theory of Banach spaces and plays a key role in the development of the Banach space theory. Our aim is to investigate some new results on complemented subspaces, to present a history of the subject, and to introduce some open problems. 
\end{abstract}
\newpage

\begin{center}
\section {Introduction.}
\end{center}

The problem related to complemented subspaces are in the heart of the theory of Banach spaces. These are more than fifty years old and play a key role in the development of the Banach space theory. Our aim is to review of results on
complemented subspaces, to present a history of the subject, and to introduce
some open problems.

We start with simple observations concerning definition and properties of complemented subspaces. Some useful sources are $\cite{CON}, \cite{HEL}, \cite{M-V}$.

Let $X$ be a normed space, $M$, $N$ be algebraically complemented subspaces of $X$ (i.e. $M+N=X$ and $M \cap N=\{0\}$), $\pi:X \to \frac{X}{M}$ be the quotient map,
$\phi:M \times N \to X$ be the natural isomorphism $(x,y) \mapsto x+y$ and $P:X \to M, P(x+y)=x , x \in M, y \in N$ be the projection of $X$ on $M$
along $N$. Then the following statements are equivalent:

(i) $\phi$ is a homeomorphism.

(ii) $M$ and $N$ are closed in $X$ and $\pi|_N$ is a homeomorphism.

(iii) $M$ and $N$ are closed and $P:X \to M$ is a bounded projection.

The Subspaces $M$ and $N$ are called topologically complemented or simply complemented if each of the above equivalent statements holds. If $N_1, N_2$ are complemented subspaces of a closed subspace $M$, then $N_1$ and $N_2$ are isomorphic Banach spaces. \\
It is known that every finite dimensional subspace is complemented and every algebraic complement of a finite codimension subspace is topologically complemented.

In a Banach space $X$, applying the closed graph theorem we can establish that two closed subspace are algebraically complemented if and only if they are complemented. Moreover, if $M$ is a closed subspace of $X$, then $M$ is complemented if and only if the following equivalent assertions hold:

(I) The quotient map $i:M \hookrightarrow X$ has a left inverse as a continuous operator .

(II) The natural projection $\pi:M \to \frac{X}{M}$ has a right inverse as a continuous operator.

$l^{\infty}$ is complementary in every normed space $X$ containing it isomorphically as a closed subspace $\cite{M-V}$. Also, If $c_\circ$ is subspace of a separable Banach space $X$ , then there is a bounded projection $P$ of $X$ onto $c_\circ$ of norm $\leq 2$, cf. $\cite{SOB}$.

Suppose now that $F$ is a retract of a Banach space $X$, i.e. $F$ is a Banach subspace of $X$ and there is a continuous linear map $\phi :X \to F$ such that for all $x \in F, \phi(x)=x$. Then $C_\circ(X-F)=\{f \in C(X): f(x)=0 {\rm ~for~ all~} x\in F\}$ is complemented in $C(X)$. In fact, by defining $P:C(X) \to C(X)$ by $P(g)=g\circ \phi$, we have $P^2=P, \|P(g)\|= _{\stackrel{Sup}{x \in X}}|g(\phi(x))|\leq\|g\|$ and $Ker P=\{ g \in C(X) | g(\phi(x))=0 {\rm~ for~ all~} x\in X\}=C_\circ(X-F)$.\\
Hence we may say that "complemented ideal' is the Gelfand dual of "retract closed subspace" (see $\cite{MOS}$).

There are non-complemented closed subspaces. For example, let $X$ be the disk algebra, i.e. the space of all analytic functions on $\{z\in {\bf C}; |z|<1\}$ which are continuous on the closure of $D$. Then the subspace of $C(T)$ consisting of the restictions of functions of $X$ to $T=\{z\in{\bf C}; |z|=1\}$ is not complemented in $X$ (see $\cite{HOF}$).

Throughout the paper $c_\circ, c, l_\infty, l_p$ denot the space of all complex sequences $\{x_n\}$ such that $\displaystyle{\lim_{n\to \infty}}x_n=0,~ \{x_n\}$ is convergent, $\{x_n\}$ is bounded, and $\displaystyle{\sum_{n=1}^\infty}|x_n|^p <\infty$, respectively. In addition, $L_p$ denotes the $L_p$-space over the Lebesgue interval $[0,1]$. The reader is referred to $\cite{J-L2}$ and $\cite{L-T}$ for undefined terms and notation.

\begin{center}
\section {Complementary subspace problem and related results.}
\end{center}

This problem asks, in general, which closed subspaces of a Banach space are complemented?

In 1937, Murray $\cite{MUR1}$ proved, for the first time, that $l_p, p\ne 2, p>1$ has non-complemented subspace.

Phillips $\cite{PHI}$ proved that $c_\circ$ is non-complemented in $l^{\infty}$. This significant fact has been refined, reproved or generalized by many mathematicians, cf. $\cite{PEL2}, \cite{GRO}, \cite{ROS2}$ and $\cite{N-K}$.

Banach and Mazur showed that all subspaces in $C[0,1]$ which are isometrically isomorphic to $l_1$ or $L^1[0,1]$ are non-complemented, cf. $\cite{SEM}$ and $\cite{BAN}$.

In 1960, Pelczynski $\cite{PEL1}$ showed that complemented subspaces of $l_1$ are isomorphic to $l_1$. K\"othe $\cite{KOT}$ generalized this result to the non-separable case.

In 1967, Lindenstrauss $\cite{LIN2}$ proved that every infinite dimensional complemented subspace of $l^{\infty}$ is isomorphic to $l^{\infty}$. This also holds if $l^{\infty}$ is replaced by $l_p$, $1 \leq p<\infty$, $c_\circ$ or $c$.

It is shown by Lindenstrauss $\cite{LIN1}$ that if the Banach space $X$ and its closed subspace $Y$ are generated by weakly compact sets (in particular, if $X$ is reflexive), then $Y$ is complemented in $X$.

In 1971, Lindenstrauss and Tzafriri $\cite{L-T}$ proved that every infinite dimensional Banach space which is not isomorphic to a Hilbert space contains a closed non-complemented subspace.

Johnson and Lindenstrauss $\cite{J-L1}$ proved the existence of a continuum of non-isomorphic separable ${\cal L}^1$-spaces. (An ${\cal L}^1$-space is a space $X$ for which $X^{**}$ is a complemented subspace of an $L^1$-space)

Classically known complemented subspaces of $L_p, 1<p<\infty ~,p\ne 2$ are $l_p, l_2, l_p\oplus  l_2$ and $L_p$ itself. In 1981, Bourgain, Rosenthal and Schechtman $\cite{B-R-S}$ proved that up to isomorphism, there exist uncountably many complemented subspaces of $L_p$.

It is shown that a complemented subspace $M$ of $l_{\infty}^*$ is isomorphic to $l_{\infty}^*$ provided $M$ is either $w^*$-closed or isomorphic to a bidual space, cf. $\cite{MET}$.

Pisier $\cite{PIS}$ established that any complemented reflexive subspace of a $C^*$-algebra is necessarily linearly isomorphic to a Hilbert space.

In 1993, Gowers and Maurey $\cite{G-M1}$ showed that there exists a Banach space $X$ without non-trivial complemented subspaces.

If $E$ is one of the spaces $l_p,~(1\leq p\leq \infty)$ or $c_\circ$, and $X$ is a vector space complemented in $E$ which contains a vector subspace $Y$ complemented in $X$ and isomorphic to $E$, then $X$ is isomorphic to $E$. Moreover, each infinite dimensional vector subspace complemented in $E$ is isomorphic to $E$. Conversely, if $Y$ is a vector subspace of $E=l^2$ or $c_\circ$ which is isomorphic to $E$, then $Y$ is complemented in $E$.

If $X$ is an infinite dimensional vector subspace complemented in some space $C(S)$, then $X$ contains a vector subspace isomorphic to $c_\circ$.

Randrianantoanina $\cite{RAN}$ showed that if $X$ and $Y$ are isometric subspaces of $L_p ~(p\ne 4,6,...)$, and $X$ is complemented in $L_p$ then so is $Y$. Moreover, the projection constant does not change. This number is defined to be $\inf\{\|T\|  : T:L_p \to X$ is a bounded linear projection of $L_p$ onto $X\}$.

The above theorem fails in the case $p\geq 4$ is an even integer, i.e. there exist pairs of isomorphic subspaces $X$ and $Y$ of $L_p$ to itself so that $X$ is complemented and $Y$ is not.

\begin{center}
\section {Schroeder-Bernstein Problem.}
\end{center}

If two spaces are isomorphic to complemented subspaces of each other, are then they isomorphic?

There are negative solutions to this problem.(see $\cite{GOW}$ and $\cite{G-M}$)

\begin{center}
\section {Basis and complemented subspaces.}
\end{center}

A Schauder basis for a Banach space $X$ is a sequence $\{x_n\}$ in $X$ with the property that every $x \in X$ has a unique representation of the form $x=\displaystyle{\sum_{n=1}^{\infty}}\alpha_nx_n ; \alpha_n \in {\bf C}$ in which the sum is convergent in the norm topology, cf. $\cite{J-L2}$. For example, the trigonometrical system is a basis in each space $L^p[0,1] , 1<p<\infty$.

Pelczynski $\cite{PEL1}$ showed that any Banach space with a basis is a complemented subspace of an isomorphically unique space.

In 1987, Szarek $\cite{SZA}$ showed that there is a complemented subspace without basis of a space with a basis and answered therefore to a problem of fifty years old.

\begin{center}
\section {Approximation property and complemented subspaces.}
\end{center}

A Banach space $X$ has the approximation property (AP) if for every $\epsilon>0$ and each compact subset $K$ of $X$ there is a finite rank operator $T$ in $X$ such that for each $x \in K , \|Tx-x\|<\epsilon$. If there is a constant $C>0$ such that for each such $T$, $\|T\|\leq C$, then $X$ is said to have bounded approximation property (BAP), cf. $\cite{J-L2}$. For example, every Banach space with a basis has BAP.

Pelczynski $\cite{PEL1}$ proved that every Banach space with the BAP can be complementably embedded in a Banach space with a basis.

\begin{center}
\section {Complemented minimal subspaces.}
\end{center}

A Banach space $X$ is called minimal if every infinite dimensional subspace $Y$ of $X$ contains a subspace $Z$ isomorphic to $X$. For example $c_\circ$ is minimal. If $Z$ is also complemented then $X$ is said to be complementary minimal. Casazza and Odell $\cite{C-O}$ showed that Tsirelson's space $T$ (see $\cite{TSI}$ and $\cite{F-J}$) have no minimal subspaces.

Casazza, Johnson and Tzafriri $\cite{C-J-T}$ showed that the dual $T^*$ of $T$ is minimal but not complementary minimal.

\begin{center}
\section {quasi-complemented subspaces.}
\end{center}

A closed subspace $Y$ of a Banach space $X$ is called quasi-complemented if there exists a closed subspace $Z$ of $X$ such that $Y \cap Z=\{0\}$ and $Y+Z$ is dense in $X$.

Then such a subspace $Z$ is said to be a quasi-complement of $Y$. Those notions are first introduced by Murray $\cite{MUR2}$.

Every closed subspace of $l_\infty$ is quasi-complemented, cf. $\cite{ROS2}$. Also Mackey $\cite{MAC}$ proved that in a separable Banach space every subspace is quasi-complemented.

Rosenthal $\cite{ROS1}$ showed that if $X$ is a Banach space, $Y$ is a closed subspace of $X$, $Y^*$ is $W^*$-separable and the annihilator $Y^{\perp}$ of $Y$ in $X^*$ has an infinite dimensional reflexive subspace, then $Y$ is quasi-complement in $X$.

\begin{center}
\section {Weakly complemented subspaces.}
\end{center}

A closed subspace of a Banach space $X$ is called weakly complemented if the dual $i^*$ of the natural embedding $i:M \hookrightarrow X$ has a right inverse as a bounded operator.

For example, $c_\circ$ is weakly complemented in $l_\infty$, not complemented in $l_\infty$ (see $\cite{WHI}$).

If $M$ is complemented in $X$ with the corresponding projection $P$, then the adjoint of $id_X-P$ is a projection in $B(X)$ with the
range $M^o=\{ f \in X^*; f|_M=0 \}$. Hence $M$ is weakly complemented in $X$.

\begin{center}
\section {contractively complemented subspaces.}
\end{center}

As mentioned before, a closed subspace $Y$ of a Banach space $X$ is said to be complemented if it is the range of a bounded linear projection $P:X \to X$. If $\|P\|=1$, $Y$ is called a contractively complemented or $1$-complemented subspace of $X$.

Let $X$ be a Banach space with ${\rm dim} X \geq 3$. Then $X$ is isometrically isomorphic to a Hilbert space iff every subspace of $X$ is the range of a projection of norm 1 (see $\cite{KAK}$ and $\cite{BOH}$).

In 1969, Zippin $\cite{ZIP}$ proved that every saparable infinite dimensional $L_1$-predual space (i.e a Banach space whose dual is isometric to $L_1(\mu)$ for some measure space $(\Omega,\Sigma,\mu)$ )) contains a contractively complemented subspace isomorphic to $c_\circ$.

Lindenstrauss and Lazar $\cite{L-L}$ proved that $X$ contains a contractively complemented subspace isometric to some space $C(S)$ when $X^*$ is non-separable.

{\bf Question.} Let X be a Banach space and $T:X \to X$ be an isometry. Is the range of $T$ is contractively complemented in $X$?\\
In Hilbert and $L^p , (1\leq p < \infty)$ spaces, we have an affirmative answer. In case $C[0,1]$, however, it may happen that the range of an isometry is not complemented, cf. $\cite{DIT}$.

Pisier $\cite{PIS}$ proved that if $M$ is a Von Neumann subalgebra of $B(H)$ which is complemented in $B(H)$ and isomorphic to $M \otimes M$, then $M$ is contractively complemented.

\begin{center}
\section {Prime Banach spaces and complemented subspaces.}
\end{center}

A Banach space $X$ is called prime if each infinite dimensional complemented subspace of $X$ is isomorphic to $X$, cf. $\cite{L-T}$.

Pelczynski $\cite{PEL1}$ proved that $c_\circ$ and $l_p~ (1 \leq p < \infty)$ are prime. Lindenstrauss $\cite{LIN2}$ proved that $l^{\infty}$ is also prime. Gowers and Maurey $\cite{G-M1}$ constructed some new prime spaces.

\begin{center}
\section {Complemented subspaces of topological products  and sums of Banach spaces.}
\end{center}

Metafune and Moscatelli $\cite{M-M}$ proved that when $X$ is one of the Banach spaces $l_p  (1\leq p\leq\infty)$ or $c_\circ$, then each infinite dimensional complemented subspace of $X^N$ is isomorphic to one of the spaces $\omega , \omega\times X^N$ or $X^N$, where $\omega=K^N$ (K is the scaler field) and $X^N$ is the product of countably many copies of $X$.

In $\cite{D-O}$, the authors obtained a complete description of the complemented subspace of the topological product $l_{\infty}^m$ where m is an arbitrary cardinal number.

Every complemented subspace of a product $V=\prod_{i \in I}X_i$ of Hilbert spaces is isomorphic to a product of Hilbert spaces (I is a set of arbitrary cardinal), cf. $\cite{DOM}$.

Ostraskii $\cite{OST}$ showed that not all complemented subspaces of countable topological products of Banach spaces are isomorphic to topological products of Banach spaces.

Chigogidze $\cite{CHI1}$ proved that complemented subspaces of a locally convex direct sum of arbitrary collection of Banach spaces are isomorphic to locally convex direct sum of complemented subspaces of countable subsums.

Chigogidze $\cite{CHI2}$ proved that a complemented subspace of an uncountable topological product of Banach spaces is isomorphic to a topological product of complemented subspaces of countable subproducts and hence isomorphic to a topological product of Frechet spaces.

\begin{center}
\section {Some interesting problems.}
\end{center}

The following problems in this area arise:\\

1) Given a Banach space $X$, characterize the isomorphic types of its complemented subspaces.

2) Given a Banach space $X$, characterize the isomorphic types of such Banach space $Z$ that every vector subspace of $Z$ isomorphic to $X$ is complemented in $Z$.

3) Is every complemented vector subspace of $C(S)$ isomorphic to some $C(S_1)$?

4) If a Banach space $X$ is complemented in every Banach space containing it, is $X$ isomorphism to some $C(S)$ over a Stone space S? (A space is Stonian if the closure of every open set is open)

5) Does every complemented subspace of a space with an unconditional basis have an unconditional basis? Recall that an unconditional basis for a Banach space is a basis $\{x_n\}$ such that every permutation of $\{x_n\}$ is also a basis or equivalently, the convergence of $\sum \alpha_n x_n$ implies the convergence of every rearrangement of the series, cf. $\cite{J-L2}$.

6) If a von Neumann algebra is a complemented subspace of $B(H)$, is it then injective?

7) Are $l_p, 1 \leq p \leq \infty$ and $c_\circ$ the only prime Banach spaces with an unconditional basis? is still open.

{\bf Remark.} Some pieces of information are taken from Internet-based resources without mentioning the URL's.

\end{document}